\documentclass[12pt]{article}
\usepackage{amsmath}
\usepackage{amsfonts}
\usepackage{mathrsfs}
\usepackage{amssymb}
\def\BC{{\tilde{\cal C}}}

\def\Vir{\mbox{Vir}}

\def\a{\alpha}
\def\b{\beta}
\def\d{\delta}
\def\D{\Delta}

\def\L{\Lambda}

\def\sc{\scriptstyle}
\def\ssc{\scriptscriptstyle}

\def\cl{\centerline}
\def\ol{\overline}

\def\bs{\backslash}

\def\vs{\vspace*}
\def\ni{\noindent}

\def\WW{{\cal W}}
\def\BB{{\cal B}}
\def\Z{\mathbb{Z}}
\def\sZ{{\ssc\,}\mathbb{Z}}
\def\C{\mathbb{C}}

\numberwithin{equation}{section}
\newtheorem{theo}{Theorem}[section]
\newtheorem{defi}[theo]{Definition}

\newtheorem{coro}[theo]{Corollary}
\newtheorem{lemm}[theo]{Lemma}

\usepackage[usenames]{color}

\begin{document}

\cl{\Large\bf {Classification of quasifinite representations}}
\cl{\Large\bf { of a Lie algebra related to Block type
}\footnote{\indent \ \ Supported by NSF grant 10825101 of China, the Fundamental Research Funds for the Central Universities, Shanghai Municipal Science and Technology Commission grant (12XD1405000)}} \vskip5pt

\cl{{Yucai Su$^{\,*}$, Chunguang Xia$^{\,\dag}$, Ying Xu$^{\,\ddag}$}}

\cl{\small\it $^{\,*}$Department of Mathematics, Tongji University, Shanghai 200092, China}

\cl{\small\it $^{\,\dag}$Department of Mathematics, China University of Mining and \vs{-4pt}Technology}
\cl{\small\it Xuzhou 221116, Jiangsu, China}

\cl{\small\it $^{\,\ddag}$Department of Mathematics, Hefei University of Technology, Hefei 230009, Anhui,  China}

\cl{\small\it Email: ycsu@tongji.edu.cn, chgxia@mail.ustc.edu.cn,
xying@mail.ustc.edu.cn} \vs{5pt}\par\ni {\small{\bf Abstract.} A
well-known theorem of Mathieu's states that a Harish-chandra
module over the Virasoro algebra is either a highest weight
module, a lowest weight module or a module of the intermediate
series. It is proved in this paper that an analogous result also
holds for the Lie algebra $\BB$ related to Block type, with basis
$\{L_{\a,i},C\,|\,\a,i\in\Z, i\ge0\}$ and relations
$[L_{\a,i},L_{\b,j}]=((i+1)\b-(j+1)\a)L_{\a+\b,i+j}
+\d_{\a+\b,0}\d_{i+j,0}\frac{\a^3-\a}{6}C$. %, $[C,L_{\a,i}]=0$.
Namely, an irreducible quasifinite $\BB$-module is either a
highest weight module, a lowest weight module or a module of the
intermediate series.\\[4pt]
} \vskip10pt

\ni {\bf 1. \
Introduction}\setcounter{section}{1}\setcounter{theo}{0}
%\section{Introduction}

\ni Since a class of infinite dimensional simple Lie algebras was
introduced by Block \cite{B}, generalizations of Lie algebras of this
type (usually referred to as {\it Lie algebras of Block type}) have
been studied by many authors (see for example,
%\cite{DZ, OZ, S4, S5, WT1, WT2, X1, X2, ZZ, Z, ZM}
[\ref{DZ}, \ref{OZ}, \ref{S4}--\ref{X2}, \ref{ZZ}--\ref{ZM}]). Lie algebras of Block type are closely
related to the Virasoro algebra, Virasoro-like algebra, or special
cases of (generalized) Cartan type $S$ Lie algebra or Cartan type
$H$ Lie algebra (e.g., \cite{X3}). It is well known that although Cartan
type Lie algebras have a long history, their representation theory
is however far from being well developed. In order to better
understand the representation theory of Cartan type Lie algebras, it
is very natural to first study representations of special cases of
Cartan type Lie algebras. Partially due to these facts, the study of
Lie algebras of this kind has recently attracted some authors'
attentions.

The author in [\ref{S4}--\ref{S6}] %\cite{S4, S5, S6}
presented a classification of the so-called
{\it quasifinite modules} (which are simply $\Z$-graded modules with
finite dimensional homogenous subspaces) over some Block type Lie
algebras. In particular, he  studied the representations of the
Block type Lie algebra $\ol{\BB}$ with basis
$\{L_{\a,i},C|\,\a,i\in\Z,i\geq-1\}$ over $\C$ and relations \cite{S4}
\begin{equation*}%\label{bar-B-block}
[L_{\a,i},L_{\b,j}]=((i+1)\b-(j+1)\a)L_{\a+\b,i+j}+
\a\d_{\a+\b,0}\d_{i+j,-2}C, \ \ [C,L_{\a,i}]=0,
\end{equation*}
for $\a,\b\in\Z,i,j\ge-1$. It is shown that quasifinite modules
with dimensions of homogenous subspaces being uniformly bounded
are all trivial --- this renders the representations of this kind
do not seem to be much interesting. On the other hand, the authors
in \cite{WS} considered Verma type modules of the above Lie algebra
$\ol{\BB}$. However it turns out that these Verma type modules,
regarded as $\Z$-graded modules, are all with infinite dimensional
homogenous subspaces (except the subspace spanned by the generator
which has dimension $1$). Probably because of the above, the
authors in \cite{WT1} turned to study representations of the Lie
algebra $\BB$ related to Block type, which is a subalgebra of
$\ol{\BB}$, with basis $\{L_{\a,i},C\,|\,\a,i\in\Z,i\ge0\}$ and
relations
\begin{equation}\label{B-block}
[L_{\a,i},L_{\b,j}]=((i+1)\b-(j+1)\a)L_{\a+\b,i+j}
+\d_{\a+\b,0}\d_{i+j,0}\frac{\a^3-\a}{6}C.
\end{equation}
This Lie algebra $\BB$, as stated in \cite{WT1}, is interesting in the
sense that it contains the following subalgebra
\begin{equation}\label{Viraaa}\Vir={\rm span}\{L_{\a,0},C\,|\,\a\in\Z\},\end{equation} which is isomorphic to the
well-known Virasoro algebra (cf.~\eqref{Vir}), while the Lie algebra $\ol{\BB}$
contains no such a subalgebra. Due to this, one will see that the
representation theory for the Lie algebra  $\BB$ is quite
different from that for $\ol{\BB}$. Moreover, $\BB$ is related to
the well-known $W$-infinity Lie algebra $\WW_\infty$ (e.g., \cite{SX}) in the
following way: Recall that the $W$-infinity Lie algebra
$\WW_{1+\infty}$ is defined to be the universal central extension of
infinite dimensional Lie algebra of differential operators on the
circle, which has a basis $\{ x^\a D^i, C|\,\a,i\in\Z,i\geq0\}$
with $D=\frac{d}{dx}$, and relations
$$
[x^\a D^i,x^\b D^j]=x^{\a+\b}((D+\b)^iD^j-D^i(D+\a)^j)+
\d_{\a+\b,0}(-1)^ii!j! \binom{\a+i}{i+j+1}C.
$$
Then the $W$-infinity algebra $\WW_\infty$, the universal central
extension of infinite dimensional Lie algebra of differential
operators on the circle of degree at least one, is simply the
subalgebra of $\WW_{1+\infty}$ spanned by $\{ x^\a
D^i,C\,|\,\a,i\in\Z,i\geq1\}$. If we define a natural filtration of
$\WW_\infty$ by
$$
\{0\}=(\WW_\infty)_{[-2]}\subset(\WW_\infty)_{[-1]}\subset\cdots\subset
\WW_\infty,
$$
where $(\WW_\infty)_{[-1]}=\C C$, and $%\begin{eqnarray*}
(\WW_\infty)_{[n]} = {\rm span}\{x^\a D^i,C\,|\,\a\in\Z, ~i\le
n+1\}\mbox{ for }n\ge0, $ %\end{eqnarray*}
then $\BB$ is simply the
associated graded Lie algebra of the filtered Lie algebra
$\WW_\infty$. As stated in \cite{S1, S3, SX}, the $W$-infinity algebras
arise naturally in various physical theories, such as conformal
field theory, the theory of the quantum Hall effect, etc.; among
them the $\WW_{\infty}$ algebra and $\WW_{1+\infty}$ algebra, of
interest to both mathematicians and physicists, have received
intensive studies in the literature.

We find the Lie algebra $\BB$ is interesting in another aspect that it also contains many (finitely) $\Z$-graded subquotient
algebras $\tilde\BB_{m,n}$ for $n{\sc\!}\ge{\sc\!} m{\sc\!}\ge{\sc\!}0$, where
\begin{equation}\label{bb-m}\tilde\BB_{m,n}=\BB_m/\BB_{n+1},\ \ \ \BB_m={\rm span}\{L_{a,i}\,|\,\a\in\Z,\,i\ge m\}.\end{equation} For instance, $\tilde\BB_{0,0}$ is the Virasoro algebra (thus the Virasoro algebra is both a subalgebra and a quotient algebra of $\BB$), $\tilde\BB_{0,1}$ is the $W$-algebra $W(-2,0)$, and $\tilde\BB_{1,2}$ is the Lie algebra of invariance of the free
Schr\"oinger equation (and thus $\tilde\BB_{0,2}$ is closely related to the twisted Schr\"oinger-Virasoro algebra).
It is well-known that the category of quasifinite modules (cf.~Definition \ref{defi-module}) over a quotient Lie algebra $\tilde\BB_{0,n}$
is a full subcategory of the category of quasifinite $\BB$-modules. Thus a classification of irreducible quasifinite
$\BB$-modules also gives a classification of irreducible quasifinite
$\tilde\BB_{0,n}$-modules for all $n>0$.

The authors in \cite{WT1} prove that an irreducible quasifinite
$\BB$-module is either a highest weight module,
a lowest weight
module or a uniformly bounded module. They also obtain a complete description of irreducible
quasifinite highest weight modules and modules of the intermediate series
over $\BB$. Motivated by a well-known result of Mathieu's in \cite{M}, it
is very natural to consider the classification of irreducible
quasifinite $\BB$-modules. In this paper, by proving that an
irreducible quasifinite uniformly bounded $\BB$-module is a module
of the intermediate series, we obtain the following main result.

\begin{theo}\label{MainTheo}Any irreducible quasifinite
$\BB$-module is either a highest weight module, a lowest weight
module or a module of the intermediate series.\end{theo}

As a by-product of our proof, we also obtain the following, which recovers a result in \cite[Theorem 4.8]{WT1}.
\begin{coro}\label{Maincoro}
A module of the intermediate series over $\BB$ is simply a module of the intermediate series
over the Virasoro algebra $\Vir$ with the trivial action of $\BB_1$.
\end{coro}

The proof of \cite[Theorem 4.8]{WT1} given in \cite{WT1} involves some heavy
computations. In our case here, some new techniques are employed, which renders the proof seems to be more elegant, as one may see in Section 3.
We would like to point out that the techniques used here may be used to dealing analogous problems of Lie algebras which are closely related to the Lie algebra $\BB$.
This is also our main motivation to present this paper.

Since lowest weight modules are duals of highest weight modules,
Theorem \ref{MainTheo} together with Corollary \ref{Maincoro} gives a complete
classification of irreducible quasifinite $\BB$-module. The
analogous results to the this theorem for the Virasoro algebra,
$W$-infinity algebras, higher rank Virasoro algebras, and some Lie
algebras of Block type were obtained in
[\ref{LZ}, \ref{M}, \ref{S2}--\ref{SX}]. %\cite{LZ,M,S2,S3,S4,S5,SX}.
 \vskip15pt

\ni {\bf 2. \
Preliminaries}\setcounter{section}{2}\setcounter{theo}{0}\setcounter{equation}{0}
%\section{Preliminaries}

\ni It is well-known that the  Virasoro algebra $\Vir$ is the Lie
algebra with basis $\{L_i,C\,|\,i\in\Z\}$ satisfying the relations
\begin{equation}\label{Vir}
[L_i,L_j]=(j-i)L_{i+j}+\frac{i^{3}-i}{12}\delta_{i,-j}C,\ \ \ \
[L_i,C]=0\mbox{\ \ \ for \ \ }i,j\in \Z.
\end{equation}
 A $\Vir$-module of the intermediate series must be one of $A_{a,b},
A(a), B(a), a, b\in\C$, or one of the quotient submodules, where
$A_{a,b}, A(a), B(a)$ all have a basis $\{x_k\,|\,k\in\Z\}$ such
that $C$ acts trivially and
\begin{eqnarray}\label{condition}
\!\!\!\!\!\!\!\!\!\!\!\!&\!\!\!\!\!\!\!\!\!\!\!\!\!\!\!& A_{a,b}:\
L_ix_k=(a+k+bi)x_{i+k},
\\\label{condition+}
\!\!\!\!\!\!\!\!\!\!\!\!&\!\!\!\!\!\!\!\!\!\!\!\!\!\!\!&
A(a): L_ix_k=(i+k)x_{i+k}\ \ (k\neq0), \ \ \ \ \ \ \ \
L_ix_0=i(i+a)x_i,
\\\label{condition++}
\!\!\!\!\!\!\!\!\!\!\!\!&\!\!\!\!\!\!\!\!\!\!\!\!\!\!\!&
B(a): L_ix_k=kx_{i+k}\ \ (k\neq-i),
\ \ \ \ \ \ \ \ \ \ \ \ \ L_ix_{-i}=-i(i+a)x_0,%\nonumber
\end{eqnarray}
for $i,k\in\Z.$ One also has
\begin{eqnarray*}&\!\!\!\!\!\!\!\!\!\!\!\!\!\!\!\!\!\!\!\!\!\!\!\!& A_{a,b} {\rm \ is\  irreducible}\ \
\Longleftrightarrow\ \  a \notin\Z {\rm \ or\ } a \in\Z,\ b\neq0,1,
\nonumber
\\ \label{isomo}
&\!\!\!\!\!\!\!\!\!\!\!\!\!\!\!\!\!\!\!\!\!\!\!\!&  A_{a,1}\cong
A_{a,0}\mbox{ if }a\notin\Z,\\ \label{isomo1}
&\!\!\!\!\!\!\!\!\!\!\!\!\!\!\!\!\!\!\!\!\!\!\!\!& A_{a,b}\mbox{ is
reducible, and }A'_{a,1}\cong A'_{a,0}\cong A'_{0,0}\mbox{ if
}a\in\Z,\ b=0,1,
\end{eqnarray*}
where in general, $A'_{a,b}$ denotes the unique nontrivial
composition factor of $A_{a,b}$.

Now consider the Lie algebra $\BB$. We can realize it in the space
$\C[x,x^{-1}]\otimes t\C[t]\oplus\C C$ with the bracket
\begin{equation}\label{B-block2}
[x^\a f(t),x^\b g(t)]=x^{\a+\b}(\b f'(t)g(t)-\a f(t)g'(t))
+\d_{\a+\b,0}\frac{\a^3-\a}{6}{\rm Res}_tt^{-3}f(t)g(t)C,
\end{equation}
for $\a,\b\in\Z,f(t),g(t)\in t\C[t]$, where the prime stands for
${d\over dt}$, and ${\rm Res}_tf(t)$ stands for the residue of the
Laurent polynomial $f(t)$, namely the coefficient of $t^{-1}$ in
$f(t)$. The Lie algebra  $\BB$ has a natural $\Z$-gradation
$\BB=\oplus_{\a\in\sZ}\BB_\a$ with $ \BB_\a=\{x^\a f(t)\,|\,f(t)\in
t\C[t]\}+\d_{\a,0}\C C. $ Putting $\BB_{\pm}
=\oplus_{\pm\a>0}\BB_\a$, we have the  triangular decomposition $
\BB=\BB_-\oplus\BB_0\oplus\BB_+. $ Note that $\BB_0=t\C[t]\oplus\C
C$ is an infinite dimensional commutative subalgebra of $\BB$ (but
not a Cartan subalgebra). When we study representations of a Lie
algebra of this kind, as pointed in [\ref{S4}--\ref{S6}], %\cite{S4,S5},
we encounter the
difficulty that though it is $\Z$-graded, the graded subspaces are
still infinite dimensional, thus the study of quasifinite modules is
a nontrivial problem.

Denote
\begin{equation}\label{denote-L==}
L_\a=x^\a, \ \ \ L_{\a,i}=x^\a
t^{i+1}\;\;\;\;\mbox{for}\;\;\;\a,i\in\Z,\,i\ge0.
\end{equation}
Then (\ref{B-block2}) is equivalent to (\ref{B-block}) and ${\rm span}\{L_\a,C\,|\,\a\in\Z\}\cong\Vir$ (cf.~\eqref{Viraaa} and \eqref{Vir}).
%In
%particular, we observe the following simple but rather deep fact
%which will be frequently used in the proof of the main theorem in a
%crucial way,
%\begin{equation}\label{deep}[L_{-1},L_{-2,1}]=0.
%\end{equation}
%Denote $ %\begin{equation}\label{denote-Vir==}
%\Vir={\rm span} \{L_\a,C\,|\,\a\in\Z\},
%%\mbox{ \ where }L_\a=L_{\a,0}
% %\end{equation}
%$ which forms a Virasoro subalgebra under bracket (\ref{B-block})
%(cf.~\eqref{Vir}).

\begin{defi}\label{defi-module}\rm
A module $V$ over $\BB$ is called
\begin{itemize}\parskip-3pt
\item {\it $\Z$-graded} if $V=\oplus_{\a\in\Z}V_\a$ and $\BB_\a
  V_\b\subset V_{\a+\b}$ for all $\a,\b$;
\item {\it quasifinite} if
  it is $\Z$-graded and ${\rm dim\ssc\,}V_\b\! <\!\infty$ for all $\b$;
\item {\it uniformly bounded} if it is $\Z$-graded and there is $N\!>\!0$
  such that ${\rm dim\ssc\,}V_\b\!\le\! N$ for all $\b$;
\item a {\it module of the
  intermediate series} if  it is $\Z$-graded and ${\rm
  dim\ssc\,}V_\b\le 1$ for all $\b$;
\item a {\it highest}
  (respectively~{\it lowest}) {\it weight module} if there exists some
  $\L\in\BB_0^*$ (the dual space of $\BB_0$) such that $V=V(\L)$,
  where $V(\L)$ is a module generated by a {\it highest}
  (respectively~{\it lowest}) {\it weight vector} $v_\L\in V(\L)_0$,
  i.e., $v_\L$ satisfies
$$
hv_\L=\L(h)v_\L\;\;\mbox{where}\;\;\;h\in\BB_0,\;\;\;
\mbox{and}\;\;\;\BB_+v_\L=0 \;\mbox{ (respectively $\BB_-v_\L=0$)}.
$$
\end{itemize}
\end{defi}

%A nonzero vector $v$ in a $\Z$-graded module $V$ is called {\it
%singular} or {\it primitive} if $\BB_+v=0$.
%
{\def\mu{k}Suppose now $V=\oplus_{\mu\in\Z}V_\mu$ is an
irreducible uniformly bounded $\BB$-module. For $a
%
%\z
%
\in\C$, we let
$
\mbox{$ V[
%
%\z
%
a]=\oplus_{\mu\in\Z}V_\mu[
%
%\z
%
a], \mbox{  where } V_\mu[
%
%\z
%
a]=\{v\in V_\mu\,|\,L_0v=(a+\mu
%
%\z
%
)v\}.$} $ Then obviously, $V[
%
%\z
%
a]$ as a $\BB$-submodule is a direct summand of $V$. Thus $V=V[
%
%\z
%
a]$ for some fixed $
%
%\z
%
a\in\C$, that is to say,
\begin{equation}\label{V-module}
\mbox{$V\!=\!\bigoplus\limits_{\mu\in\Z}V_\mu, \mbox{ where
}V_\mu%=V_\mu[
%
%\z
%
%a]
\!=\!\{v\in V\,|\,L_0v=(a+k)v\}\mbox{ ({\it weight space with weight $a+k$})}.$}
\end{equation}
Note that regarding as a $\Vir$-module, $V$ is also a uniformly
bounded $\Vir$-module. Therefore by the result of \cite{MP1, MP2}, we
have the following lemma.

\begin{lemm}\label{Lemm-uniformly}If $V$ is an irreducible
uniformly bounded $\BB$-module as in $(\ref{V-module})$, then there
exists a non-negative integer $N$ such that ${\rm
dim\ssc\,}V_\mu[
%
%\z
%
a]=N$ for all $\mu\in\Z$ with $\mu+
%
%\z
%
a\neq0$.
\end{lemm}

 \ni {\bf 3. \ Proof of the main
theorem}\setcounter{section}{3}\setcounter{theo}{0}\setcounter{equation}{0}
%\section{Proof of the main theorem}

\ni The aim of this section is to prove Theorem \ref{MainTheo}.
%??????????????????????
%Unfortunately, the proof has to be divided into case to case
%discussions, which will occupy the must part of the paper. The
%reason we have to make a case by case consideration can be seen from
%the results of [MP2] on a classification of indecomposable uniformly
%bounded $\Vir$-modules $V$ with ${\rm dim\ssc\,}V_k\le 2$ for all
%$k$. Since the results in [MP2] contain several different cases, and
%there is no uniform way to state the whole results,
%%which makes the result stated in [MP2] difficult to be used here,
%we find out that
%sometimes it is more convenient for us to deduce directly the
%results required by us.
%
%Suppose
Let $a$ be fixed such that \eqref{V-module} holds (we always
choose $a$ to be zero if $a\in\Z$). For any $\mu\in\Z$, the
following notation will be frequently used,
\begin{equation}\label{freq}\tilde
\mu=a+\mu.\end{equation}%
%
%
%\vskip4pt
%T \ni{\textbf{3.1. General discussions.}}
%\subsection{General discussions $aaa$}
Let $V$ as in (\ref{V-module}) be an irreducible uniformly bounded
$\BB$-module. By Lemma \ref{Lemm-uniformly}, ${\rm dim\ssc\,}V_k
%[
%
%\z
%a]
=N$ for all $
%
%\z
%
k\in\Z$ with $\tilde k
%
%\z
%
\neq0$.
Without loss of generality, we can suppose $N\geq1$.
%
%Without loss of generality, we can suppose ${\rm dim\ssc\,}V_\mu=N$
%for all $\mu\in\Z$ (otherwise, the arguments are similar to that in
%[S3, S4, S5]).
Regarding $V$ as a $\Vir$-module and choosing a composition series, by a
well-known theorem of Matheiu's \cite{M}, we can take a basis $Y_k=(y^{(1)}_k,...,y^{(N)}_k)$ of $V_k$ with $\tilde k\ne0$ satisfying
\begin{equation}\label{comp-vir}
L_{\a}Y_k=Y_{\a+k}A_{\a,k}\mbox{ \ for \ }\a,k\in\Z\mbox{ \ with \ }\tilde k,\,\a+\tilde k\ne0,
\end{equation}
such that $A_{\a,k}$ is an upper-triangular matrix with diagonals being $\tilde k+b_p\a$ for $p=1,...,N$ and some $b_p\in\C$.
Our first result is the following.
\begin{lemm}\label{fir-lemm}\begin{itemize}\parskip-3pt\item[\rm(1)]
For all $i\gg0$, the action of $L_{1,i}$ on $V$ is trivial, i.e., $L_{1,i}|_V=0.$
\item[\rm(2)]There exists some $j\ge0$ such that $\BB_{j+1}V=0$, where $\BB_{j+1}$ is defined in \eqref{bb-m}.\end{itemize}
\end{lemm}\noindent{\it Proof.~}(1) Fix $i>0$ and suppose $L_{1,i}Y_k=Y_{k+1}T_k$ for $\tilde k,\tilde k+1\ne0$ and some $N\times N$ matrix $T_k=(t^{p,q}_k)_{p,q=1}^N$ (the symbol $T_k$ remains to be undefined if $\tilde k=0,-1$). Assume $T_k\ne0$ and let $(p,q)$ be (first) the leftmost and (then) the lowermost position such that
\begin{equation}\label{t-k==0}t_k\ne0\mbox{ \  for some \ }\tilde k,\tilde k+1\ne0,\end{equation}
where in general, we denote $t_j=t^{p,q}_j$ for all possible $j$.
 Let $\a,\b\gg0$ be such  that $T_{k},T_{\a+k},$ $T_{\b+ k},T_{\a+\b+ k}$ appearing below are all defined.
Applying the equation \begin{equation}\label{equai-j}
\Big(1-(i+1)(\a+\b)\Big)[L_\a,[L_\b,L_{1,i}]]=\Big(1-(i+1)\b\Big)\Big(1+\b-(i+1)\a\Big)[L_{\a+\b},L_{1,i}],
\end{equation} to $Y_k$, we obtain
\begin{eqnarray}\label{comp-vir-L-i}
&\!\!\!\!\!\!\!\!\!\!\!\!&\Big({\sc\!}1\!-\!(i\!+\!1)(\a\!+\!\b){\sc\!}\Big)\Big({\sc\!}A_{\a,1+\b+k}(A_{\b,1+k}T_k
\!-\!T_{\b+k}A_{\b,k})\!-\!(A_{\b,1+\a+k}T_{\a+k}\!-\!T_{\a+\b+k}A_{\b,\a+k})A_{\a,k}{\sc\!}\Big)\nonumber\\
&\!\!\!\!\!\!\!\!\!\!\!\!\!\!\!\!&=\Big(1-(i+1)\b\Big)\Big(1+\b-(i+1)\a\Big)(A_{\a+\b,1+k}T_k-T_{\a+\b+k}A_{\a+\b,k}).
\end{eqnarray}
Comparing the $(p,q)$-entry, we have
\begin{eqnarray}\label{comp-vir-L-i-p,q}
&\!\!\!\!\!\!\!\!\!\!\!\!\!\!\!\!\!&\Big(1\!-\!(i\!+\!1)(\a\!+\!\b)\Big)\Big((1+\b+\tilde k+b_q\a)\big((1+\tilde k+b_q\b)t_k-t_{\b+k}(\tilde k+b_p\b)\big)\nonumber\\&\!\!\!\!\!\!\!\!\!\!\!\!\!\!\!\!\!&-\big((1+\a+\tilde k+b_q\b)t_{\a+k}-t_{\a+\b+k}(\a+\tilde k+b_p\b)\big)(\tilde k+b_p\a)\Big)\nonumber\\
&\!\!\!\!\!\!\!\!\!\!\!\!\!\!\!\!\!&=\Big(1\!-\!(i\!+\!1)\b\Big)\Big(1\!+\!\b\!-\!(i\!+\!1)\a\Big)\Big(\big(1\!+\!\tilde k\!+\!b_q(\a\!+\!\b)\big)t_k\!-\!t_{\a+\b+k}\big(\tilde k\!+\!b_p(\a\!+\!\b)\big)\Big).
\end{eqnarray}
Setting the triple $(\a,\b,k)$ in \eqref{comp-vir-L-i-p,q} to be $(\a,\a,k-\a),\,(\a,-\a,k),\,(-\a,-\a,k+\a)$ respectively, we obtain a system of three equations on variable $t_{k-\a},t_k,t_{k+\a}$.
Note that the determinant  $\D_{i,\a,k}$ of coefficients of this system
is a polynomial on $i,\a,k$ with the degree of $i$ being $\le6$. If we compute
the coefficient of $i^6$ in $\D_{i,\a,k}$, then only those terms in \eqref{comp-vir-L-i-p,q} with factor $i^2$ need to be considered, in this case, \eqref{comp-vir-L-i-p,q} can be simplified as $$0=i^2\b\a\Big(\big(1\!+\!\tilde k\!+\!b_q(\a\!+\!\b)\big)t_k\!-\!t_{\a+\b+k}\big(\tilde k\!+\!b_p(\a\!+\!\b)\big)\Big).$$ From this,
one can
easily compute that the coefficient of $i^6$ in $\D_{i,\a,k}$ is $1 + 2(\a + \tilde k) - 4\a^2( b_p +  b_p^2 +  b_q -b_q^2)$, which is nonzero when $\a\gg0$. Therefore, $\D_{i,\a,k}\ne0$
when $i,\a\gg0$, and we obtain $t_k=0$ if $i\gg0$, a contradiction with \eqref{t-k==0}. This
proves $L_{1,i}V_k=0$ for $\tilde k,\tilde k+1\ne0$ and $i\gg0$, and in particular, we have (1) if $a\notin\Z$.

Now assume $a=0$. Take a basis $Y_0=(y^{(1)}_0,...,y^{(N')}_0)$ of $V_0$ (assume ${\rm dim\,}V_0=N'$), and assume
$L_{1,i}Y_{j}=Y_0T_{j},\,j=-1,0,$ for some $N'\times N$ matrix $T_{-1}$ and some $N\times N'$ matrix $T_{0}$. Applying \eqref{equai-j} to $Y_{-\a-\b-1}$ and $Y_0$, we obtain respectively (cf.~\eqref{comp-vir-L-i})
\begin{eqnarray}\label{comp-vir-L-i+11}
T_{-1}P_{\a,\b,i}=0,\ \ Q_{\a,\b,i}T_0=0,
\end{eqnarray}
where
\begin{eqnarray*}
&\!\!\!\!\!\!\!\!\!\!\!\!&
P_{\a,\b,i}\!=\!\Big({\sc\!}1\!-\!(i\!+\!1)(\a\!+\!\b){\sc\!}\Big)A_{\b,-\b-1}A_{\a,-\a-\b-1}
\!+\!\Big({\sc\!}1\!-\!(i\!+\!1)\b{\sc\!}\Big)\Big({\sc\!}1\!+\!\b\!-\!(i\!+\!1)\a{\sc\!}\Big)A_{\a+\b,-\a-\b-1},\\
&\!\!\!\!\!\!\!\!\!\!\!\!&
Q_{\a,\b,i}\!=\!\Big({\sc\!}1\!-\!(i\!+\!1)(\a\!+\!\b){\sc\!}\Big)A_{\a,1+\b}A_{\b,1}\!-\!\Big({\sc\!}1\!
-\!(i\!+\!1)\b{\sc\!}\Big)
\Big({\sc\!}1\!+\!\b\!-\!(i\!+\!1)\a{\sc\!}\Big)A_{\a+\b,1},
\end{eqnarray*}
which are upper-triangular matrices with nonzero diagonals when $\a,\b,i\gg0$. Thus $T_0=T_{-1}=0$ by \eqref{comp-vir-L-i+11}.

(2) By (1), there exists $j$ such that $L_{1,j+1}|_V=0$. Since $\BB_{j+1}$ is an ideal of $\BB$ generated (as an ideal) by $L_{1,j+1}$, we obtain (2).
\hfill$\Box$
\begin{lemm}\label{lemm-2}
We have $\BB_1V=0$.\end{lemm}
\noindent{\it Proof.~}Let $j$ be the smallest such that Lemma \ref{fir-lemm}(2) holds. Assume $j>0$. Then
there exists some $\a\in\Z$ with \begin{equation}\label{ff====0---}L_{\a,j}|_V\ne0,\end{equation} and $V$ becomes an irreducible module over $\tilde\BB_{0,j}$ (cf.~\eqref{bb-m}). We shall use
$\tilde L_{\a,i},\,\a\in\Z$, $0\le i\le j,$ to denote basis elements of $\tilde\BB_{0,j}$ (and still use $L_\a$ to denote $\tilde L_{\a,0}$). Denote $\BC=$\linebreak ${\rm span}\{\tilde L_{\a,j}\,|\,\a\in\Z\},$ which is an ideal of $\tilde \BB_{0,j}$ and is in the center of $\tilde\BB_{1,j}$. Then \eqref{ff====0---} shows $\BC V$ is a nonzero submodule of $V$. Thus
\begin{equation}\label{CC-ff====0---}\BC V=V.\end{equation}
Take the  subspace $M\!=\!\oplus_{i=-2}^2V_i$ of $V$. Then $\tilde L_{0,j}|_M$ is a linear transformation on the finite-dimensional space $M$. Let $f(\lambda)$ be its characteristic
polynomial (or its minimal polynomial), so \begin{equation}\label{ff====0}f(\tilde L_{0,j})M=0.\end{equation} From \eqref{condition}--\eqref{condition++}, one can easily obtain by induction on $N$ that
\begin{equation}\label{V=M+VM====0}V=M+\Vir{\ssc\,}
 M,\end{equation}which can also be obtained by proving by induction on $q$ that every basis element $y^{(q)}_k\in V_k$ for all $k$ with $|k|>2$ is a linear combination of $L_{k+i}y_{-i}^{(p)}$ for $-2\le i\le 2$ and $1\le p\le q$.
 Noting that $\tilde L_{0,j}$ is in the center of the universal enveloping algebra $U(\tilde \BB_{1,j})$ of $\tilde\BB_{1,j}$,
for any polynomial $g(\lambda)$ and $\a\in\Z$, we have (cf.~\eqref{B-block})
\begin{equation}\label{ggggggggg==}
g(\tilde L_{0,j})L_\a =L_\a{\sc\,}g(\tilde L_{0,j})+[g(\tilde L_{0,j}),L_{\a}]
=L_\a{\sc\,}g(\tilde L_{0,j})+(j+1)\a g^{(1)}(\tilde L_{0,j})\tilde L_{\a,j},\end{equation}
where in general, $g^{(p)}(\lambda)$ denotes the $p$-th derivative of $g(\lambda)$. Now take $g(\lambda)=f(\lambda)^2$.
By \eqref{ff====0} and \eqref{ggggggggg==}, we have
\begin{equation}\label{ggggggggg==fff}
g(\tilde L_{0,j})L_\a M\subset 2(j+1)\a\tilde L_{\a,j}f^{(1)}(\tilde L_{0,j})f(\tilde L_{0,j})M=0\mbox{ \ for \ }\a\in\Z,\end{equation}
which together with \eqref{ff====0} and \eqref{V=M+VM====0} proves $g(\tilde L_{0,j})V=0$. From this, as in \eqref{ggggggggg==}, by considering $g(\tilde L_{0,j})L_{\a_1}\cdots L_{\a_q}$ and induction on $q$, we have
\begin{equation}\label{ggggggggg==fff-q}
g^{(q)}(\tilde L_{0,j})\tilde L_{\a_1,j}\cdots \tilde L_{\a_q,j}V=0,\end{equation}
for all $\a_1,...,\a_q\in\Z\bs\{0\}$. If, say,  $\a_1=0$, by choosing  some $\a\in\Z\bs\{-\a_2,...,-\a_q\}$, and using  $L_{0,j}=((j+1)\a)^{-1}[L_\a,L_{-\a,j}]$, we see that \eqref{ggggggggg==fff-q} still holds. From this, we see that \eqref{ggggggggg==fff-q} holds for all $\a_i\in\Z$. This in particular proves that for
$s={\rm deg\,}g(\lambda)$,
\begin{equation}\label{gg==fff-q++}
\BC^s V=0.\end{equation}
This together with \eqref{CC-ff====0---} is a contradiction. Thus $j=0$ and the result is proved.
\hfill$\Box$\vskip5pt

Lemma \ref{lemm-2} says that $V$ is simply an irreducible module over $\Vir$, thus it is a module of the intermediate series. This completes the proof of
Theorem \ref{MainTheo}.

%\vs{15pt}\par\ni {\bf References} \vs{7pt} \baselineskip=-2pt
%\parskip=-0.03truein
%\parindent=2.0cm
%\def\hang{\hangindent\parindent}
%\def\textindent#1{\indent\llap{[#1]\enspace}\ignorespaces}
%\def\re{\par\hang\textindent}
%\small\def\bf{}

\small
\footnotesize
\parskip1pt\lineskip1pt
\parskip=0pt\baselineskip=0pt
%\noindent{\bf{References}}
\def\re#1{\bibitem{#1}\label{#1}}
%\end{enumerate}
\end{document}